\documentclass{elsart}
\usepackage{amsmath,amscd,amsfonts,amssymb}

\newtheorem{Theorem}{Theorem}[section]
\newtheorem{Lemma}[Theorem]{Lemma}
\newtheorem{Corollary}[Theorem]{Corollary}
\newtheorem{Proposition}[Theorem]{Proposition}
\newtheorem{Remark}[Theorem]{Remark}
\newtheorem{Example}[Theorem]{Example}

\def\deg{\operatorname{deg}}
\def\reg{\operatorname{reg}}
\def\gin{\operatorname{in}}
\def\codim{\operatorname{codim}}

\def\sk{\smallskip}
\def\Nset{{\mathbb N}}
\def\Zset{{\mathbb Z}}
\def\Rset{{\mathbb R}}

\def\mfr{{\mathfrak m}}
\def\abf{{\mathbf  a}}
\def\bbf{{\mathbf  b}}
\def\mbf{{\mathbf  m}}
\def\nbf{{\mathbf  n}}
\def\pbf{{\mathbf  p}}
\def\qbf{{\mathbf  q}}
\def\ebf{{\mathbf  e}}
\def\xbf{{\mathbf  x}}
\def\ybf{{\mathbf  y}}
\def\tbf{{\mathbf  t}}

\def\Pcal{{\mathcal P}}
\def\Acal{{\mathcal A}}

\begin{document}

\begin{frontmatter}

\title{ Gr\"obner bases of simplicial toric ideals}
\thanks{The second  author was supported by the National Basic Research Program 
(Vietnam) and Max-Planck Institute for Mathematics in the Sciences (Germany). He
would like to thank the MIS for the financial support and the hospitality.}
\author{Michael Hellus}
\address{Universit\"at Leipzig, Fakult\"at f\"ur Mathematik und Informatik,
Augustusplatz 10/11, D-04109 Leipzig, Germany }
\ead{Michael.Hellus@math.uni-leipzig.de}
\author{L\^e Tu\^an Hoa}
\address{Institute of Mathematics Hanoi, 18 Hoang Quoc Viet Road, 10307 Hanoi, 
Vietnam}
\ead{lthoa@math.ac.vn}
\author{J\"urgen St\"uckrad}
\address{Universit\"at Leipzig, Fakult\"at f\"ur Mathematik und Informatik,
Augustusplatz 10/11, D-04109 Leipzig, Germany }
\ead{stueckrad@math.uni-leipzig.de}
\begin{abstract}  Bounds for the  maximal degree of certain Gr\"obner bases of simplicial
toric ideals are given. These bounds are close to the bound stated in Eisenbud-Goto's
Conjecture on the Castelnuovo-Mumford regularity.
\end{abstract}
 \maketitle

\end{frontmatter}

\section*{Introduction}

Let $I$ be a homogeneous ideal of a polynomial ring $R$. The coarsest measure of the 
complexity of a Gr\"obner basis (w.r.t. to a term order $\leq $) of an ideal $I$ is its maximal 
degree, which is the highest degree of a generator of the initial ideal $\gin_\leq(I)$. However, 
this quantity is not easy to be handled with. One way to study it is to use a better-behaved 
invariant, the Castelnuovo-Mumford regularity $\reg (I)$ of $I$. This invariant can be defined 
as the maximum over all $i$ of the degree minus $i$ of any minimal $i$-th syzygy of $I$, 
treating generators as 0-th syzygies. In the generic coordinates and with respect to the reverse 
lexicographic order, the degree of a minimal Gr\"obner basis of $I$ is bounded by $\reg (I)$ 
(see \cite[Corollary 2.5]{BS}). Unfortunately, this is not true for arbitrary coordinates. On the 
other hand, a famous conjecture by Eisenbud and Goto states that $\reg(I)\leq \deg (R/I) 
-\codim(R/I) + 1$, provided $I$ is a prime ideal containing no linear form (see \cite{EG}). Here 
$\deg(R/I)$ and $\codim(R/I)$ denote the multiplicity and the codimension of $R/I$, 
respectively. Thus one may guess that the {\it Eisenbud-Goto bound} $\deg (R/I) -\codim(R/I) 
+ 1$ is an expected bound for degrees of certain Gr\"obner bases of $I$.

In this paper we are interested in estimating the complexity of certain Gr\"obner bases of 
simplicial toric ideals. Toric ideals are nice, particularly because they are prime ideals 
generated by binomials. In \cite{HS} some  bounds  close to $\deg (R/I) -\codim(R/I) + 1$ 
were obtained for the Castelnuovo-Mumford regularity of simplicial toric ideals. Therefore we 
would like to pose the following problem:

 {\bf Question}: {\it Assume that $I$ is the defining  ideal of  the semigroup ring of  a 
simplicial affine semigroup $S$ over a field $K$. Does $I$ posses a Gr\"obner basis of degree 
at most $\deg K[S] -\codim K[S] + 1$ in the natural coordinates?}

A similar problem was posed by Sturmfels \cite{St1} for the class of toric ideals defined by 
so-called normal semigroups (see the question before Corollary \ref{C2}). The requirement to 
keep the natural coordinates here is essential, because then  the reduced Gr\"obner basis 
consists of binomials - which are cheep to compute and to restore. 

We will show that for many classes of simplicial  toric ideals the above question has a positive 
answer. In order to do that we first establish some general bounds for the maximal degree of 
certain Gr\"obner bases in terms of the reduction number $r(S)$ of $K[S]$ (Theorem \ref{A1}), 
or in terms of the codimension $c= \codim K[S]$ and the total degree $\alpha $ of monomials 
defining $S$ (Theorem \ref{A4}). In a lot of concrete examples these bounds are even much 
smaller then the Eisenbud-Goto bound. In the general case, we are still not able to solve the 
above problem. However, combining with a bound of \cite{HS} on $r(S)$, we can quickly show 
that $I$ possesses a Gr\"obner basis of degree at most $2(\deg K[S] -\codim K[S])$ (see 
Theorem \ref{A3}). The general bounds are given in Section \ref{A}, where examples are also 
constructed to show that  they are  close to be the best. Term orders in this section are not 
necessarily the reverse lexicographic order. In Section \ref{B} we mainly consider the reverse 
lexicographic order and derive the Eisenbud-Goto bound for certain classes of simplicial toric 
ideals. Ideals of first type come from a simple observation that degrees of their minimal 
Gr\"obner bases are bounded by the Castelnuovo-Mumford regularity if the coresponding ring  
$K[S]$ is a generalized Cohen-Macaulay ring. For ideals of second type, by using Theorem 
\ref{A4}, we can restrict ourselves to  few exceptional cases when the codimension is very big. 
Then the main technique is  to refine bounds on the reduction number or to calculate its exact 
value, so that one can apply Theorem \ref{A1}. In particular,  we show that all simplicial toric 
ideals, for which the Eisenbud-Goto conjecture is known to hold, also have a Gr\"obner basis of 
degree at most $\deg K[S] -\codim K[S] + 1$ (w.r.t. the  reverse lexicographic order).
\sk

{\it Notation}: In this paper we use bold letters to denote a vector, while their coordinates are 
written in the normal style. Thus $a_i,\ e_{1i}$ are the $i$-th coordinates of vectors $\abf,\ 
\ebf_1$, respectively; $\xbf^{\mbf} = x_1^{m_1}\cdots x_c^{m_c}$, $\ybf^{\nbf} = 
y_1^{n_1}\cdots y_d^{n_d}$ and $\tbf^{\nbf} = t_1^{n_1}\cdots t_d^{n_d}$. The  ordering 
of variables is always assumed to be $x_1 > \cdots > x_c >y_1>\cdots > y_d$. For a fixed  
term order, if not otherwise stated, then we assume that the first term in a binomial is bigger 
than the second one.

\section{Bounds}\label{A}

Let $S\subseteq \Nset^d$ be a homogeneous, simplicial affine semigroup generated by 
  a set of elements of the following type:
$$\Acal = \{\ebf_1,...,\ebf_d, \abf_1,...,\abf_c\} \subseteq M_{\alpha ,d} =\{(x_1,...,x_d)\in  
\Nset^d|\ x_1 +\cdots + x_d =\alpha \},$$
where $c\geq 2, \alpha \geq 2$ are natural numbers and $ \ebf_1 = (\alpha, 0,...,0), ..., \ebf_d 
=(0,...,0, \alpha)$.
Moreover, if $\abf_i = (a_{i1},...,a_{id})$, we can assume that the integers $a_{ij}$, where $ 
i=1,...,c ,\ j=1,...,d$, are
relatively prime.  Note that $\dim K[S] = d$ and $\codim K[S]=c$.  Let $I_\Acal$ be   the 
kernel of the homomorphism
$$\begin{array}{r}K[\xbf,\ybf] := K[x_1,...,x_c,y_1,...,y_d] \rightarrow K[S]\equiv  
K[t_1^\alpha ,...,t_d^\alpha, \tbf^{\abf_1},...,\tbf^{\abf_c}] \subseteq K[\tbf];\\ x_i \mapsto 
\tbf^{\abf_i};\ y_j \mapsto t_j^\alpha,\ i=1,...,c;\ j=1,...,d .
\end{array}$$
We call $I_\Acal$ a {\it simplicial toric ideal defined by $\Acal$} (or $S$). We will consider 
the standard grading on $K[\xbf,\ybf]$ and $K[S]$, i.e. $\deg(x_i) = \deg(y_j) =1$ and if 
$\bbf\in S$, then $\deg(\bbf) = (b_1+\cdots + b_d)/\alpha$.

Note that with respect to any term order, $I_\Acal$ has a Gr\"obner basis consisting of  
binomials (see, e.g.,  
\cite[Chapter 1]{St1}). We are interested in bounding  the maximal degree of certain 
Gr\"obner bases of $I_\Acal$. In this paper we only consider either the reverse lexicographic 
order or term orders with the following property:

 \begin{itemize} \item[(*)]  {\it If $ m,m'$ are  monomials  of $K[\xbf,\ybf] $ and $ \deg_\xbf 
m> \deg_\xbf m',$ then $ m> m',$
where $\deg_\xbf \xbf^\mbf\ybf^\nbf := m_1+\cdots + m_c$}.\end{itemize}

Let $A=A_0\oplus A_1 \oplus \cdots$, where $A_0 = K$, be a standard graded
$K$-algebra of dimension $d$. A minimal reduction of $A$ is a graded ideal $I$
generated by $d$ linear forms such that $[IA]_n = A_n$ for $n\gg 0$. The least
integer $n$ such that $[IA]_{n+1} = A_{n+1}$ is called the reduction number of
$A$ w.r.t. $I$ and will be denoted by $r_I(A)$. Note that $(t_1^\alpha,...,
t_d^\alpha)$ is a minimal reduction of $K[S]$. We denote by $r(S)$  the reduction number
 of $K[S]$ w.r.t. this minimal reduction. Then $r(S)$ is the least positive
integer $r$ such that  $(r+1)\Acal = \{e_1,...,e_d\}+r\Acal$, where for two subsets $B$
and $C$ of $\Zset^d$ we denote by $B\pm C$  the 
set of all elements of the form $b\pm c,\ b\in B, \ c\in C$, and $nB = B+\cdots + B$ ($n$ times).  
This reduction number was used in \cite{HS} to bound the Castelnuovo-Mumford regularity of 
$K[S]$. In the following result we show that one can use this number to bound the maximal 
degree of the reduced Gr\"obner basis, too.

\begin{Theorem} \label{A1} With respect to any term order specified as above, $I_\Acal$ has 
a Gr\"obner basis of degree at most $\max\{r(S)+1, 2r(S)-1\} \leq 2r(S)$.
\end{Theorem}

\begin{pf} Let $s= \max\{r(S)+1, 2r(S)-1\}$ and set
$$G =\{ \xbf^\mbf\ybf^\nbf - \xbf^\pbf\ybf^\qbf \in I_\Acal |\  \deg (\xbf^\mbf\ybf^\nbf ) = 
\deg( \xbf^\pbf\ybf^\qbf ) \leq s\}.$$
It suffices to show that $G$ is a Gr\"obner basis. Assume that this is not the case. Then one can 
find a binomial $b= \xbf^\mbf\ybf^\nbf - \xbf^\pbf\ybf^\qbf\in I_\Acal$ of the smallest 
degree $\deg b>s$ such that $\gin(g) \nmid \xbf^\mbf\ybf^\nbf $ for all $g\in G$.

If $\deg(\xbf^\mbf) \geq r(S)+1$, then we can write $\xbf^\mbf = \xbf^{\mbf'} \xbf^{\mbf"}$, 
where $\deg(\xbf^{\mbf'}) = r(S)+1$. By the definition of $r(S)$ we can find $\mbf^*, \nbf^*$ 
such that $\deg(\xbf^{\mbf^*}) =r(S)$ and  $g := \xbf^{\mbf'} - \xbf^{\mbf^*}\ybf^{\nbf^*} 
\in I_\Acal$ (note that w.r.t. both kinds of term orders,  $\xbf^{\mbf'} > 
\xbf^{\mbf^*}\ybf^{\nbf^*}$). Then $g\in G$ and $\gin(g) = \xbf^{\mbf'} \mid 
\xbf^\mbf\ybf^\nbf $, a contradiction. Thus $\deg(\xbf^\mbf) \leq r(S)$. 

If the term order satisfies the condition (*), then $\deg(\xbf^\pbf)\leq r(S)$ too. In the case of 
the reverse lexicographic order, if $\deg(\xbf^\pbf)\geq r(S)+1$, then as above, we can find 
$\pbf', \pbf"$ such that $\xbf^\pbf\ybf^\qbf - \xbf^{\pbf'}\ybf^{\pbf" + \qbf} \in I_\Acal$ and 
$\deg(\xbf^{\pbf'}) =r(S) < \xbf^{\pbf}$. Then
$$\xbf^\mbf\ybf^\nbf - \xbf^{\pbf'}\ybf^{\pbf" + \qbf} = (\xbf^\mbf\ybf^\nbf - 
\xbf^\pbf\ybf^\qbf) + ( \xbf^\pbf\ybf^\qbf - \xbf^{\pbf'}\ybf^{\pbf" + \qbf}) \in I_\Acal,$$ 
and $\xbf^\mbf\ybf^\nbf > \xbf^\pbf\ybf^\qbf > \xbf^{\pbf'}\ybf^{\pbf" + \qbf}$. Hence, 
replacing $\xbf^\pbf\ybf^\qbf $ by $\xbf^{\pbf'}\ybf^{\pbf" + \qbf}$, we may assume from the 
beginning that $\deg(\xbf^\pbf)\leq r(S)$.

Now, since $\xbf^\mbf\ybf^\nbf - \xbf^\pbf\ybf^\qbf\in I_\Acal$, we have
$$\sum_{i=1}^c m_i\abf_i + \sum_{j=1}^d n_j\ebf_j = \sum_{i=1}^c p_i\abf_i + 
\sum_{j=1}^d q_j\ebf_j.$$
From the minimality of $\deg(\xbf^\mbf\ybf^\nbf)$ we may assume that 
$\xbf^\mbf\ybf^\nbf$ and $\xbf^\pbf\ybf^\qbf$ have no common variable. That means if we 
set  $C= \{ i|\ m_i \neq 0\}$ and $D= \{j|\ n_j \neq 0\}$, then the above equality can be 
rewritten as
$$\sum_{i\in C} m_i\abf_i + \sum_{j\in D} n_j\ebf_j = \sum_{i\not\in C} p_i\abf_i + 
\sum_{j\not\in D} q_j\ebf_j.$$
Hence
$$\sum_{j\in D} \sum_{i\in C} m_ia_{ij} + \sum_{j\in D} n_j\alpha  = \sum_{j\in 
D}\sum_{i\not\in C} p_i a_{ij} = \sum_{i\not\in C} p_i \sum_{j\in D}a_{ij} \leq 
\sum_{i\not\in C} p_i\alpha .$$
This implies 
\begin{equation}\label{EA1}
\sum_{j=1}^d n_j = \sum_{j\in D} n_j \leq \sum_{i\not\in C} p_i = \deg(\xbf^\pbf).
\end{equation}
The equality holds  if  and only if $m_ia_{ij}=0$ for all $(i,j)\in C\times D$ and  
$p_ia_{ij}=0$ for all $(i,j)$ such that $ i\not\in C$ and $ j\not\in D$. This yields $\sum_{i\in 
C} m_i\abf_i = \sum_{j\not\in D} q_j\ebf_j$, which means $\xbf^\mbf - \ybf^\qbf \in 
I_\Acal$. Since $\xbf^\mbf > \ybf^\qbf$ and $\deg(\xbf^\mbf )\leq r(S)$, $g:= \xbf^\mbf - 
\ybf^\qbf \in G$. But this is impossible because $\gin(g) \mid \xbf^\mbf \ybf^\nbf $.
Hence, by (\ref{EA1}), we must have $\sum_{j=1}^d n_j < \deg(\xbf^\pbf)\leq  r(S)$, and so 
$$\deg(b) = \deg( \xbf^\mbf) + \sum_{j=1}^d n_j  \leq 2r(S)-1 \leq s,$$
a contradiction. The theorem is proved. \hfill $\square$
\end{pf}

The following example shows that  the maximal degree of a reduced Gr\"obner basis in the 
worst case must be at least $r(S)+1$.

\begin{Example} \label{A1b} {\rm Given $d\geq 2$ and $\alpha \geq d+1$. Let
\begin{equation} \label{EA1b}
\Acal = M_{\alpha ,d}\setminus \{ (\beta , \alpha -\beta ,0,...,0)|\ 2\leq \beta \leq \alpha -2\}.
\end{equation}
We denote $\abf_1 = (\alpha -1,1,0,...,0)$ and $\abf_2 = (1,\alpha -1,0,...,0)$ - the only two 
inner points of $\Acal$ in the edge $\overline{\ebf_1\ebf_2}$. If  $S \ni (\alpha -2)\abf_1 = 
\sum m_i\abf_i + \sum n_j\ebf_j$ with $\sum n_j >0$, comparing the second coordinate, one 
should have $\alpha -2 = m_1 + m_2 (\alpha -1) + n_2\alpha $. This implies $n_2=m_2 =0$ 
and 
$m_1=\alpha -2$, which is impossible, since 
$m_1= (\alpha -2) - \sum n_j <\alpha -2$. Hence $(\alpha -2)\abf_1\not\in \{\ebf_1,...,\ebf_d\} 
+ (\alpha -3)\Acal$ and $r(S) \geq \alpha -2$.

Let $\bbf = (b_1,...,b_d)\in \Nset^d$ such that $b_1 + \cdots + b_d\ \vdots \ \alpha$ and 
$b_3+\cdots +b_d>0$. By induction on $\deg(\bbf) : = (b_1+\cdots +b_d)/\alpha $, we show 
that $\bbf \in S$. The case $\deg(\bbf) =1$ follows from (\ref{EA1b}). Let $\deg(\bbf) \geq 2$. 
If $b_1\geq \alpha $, then $\bbf = \ebf_1 + \bbf'$ with $b'_3 + \cdots + b'_d >0$.  By the 
induction hypothesis, $\bbf'\in S$ and hence $\bbf\in S$. The same holds if $b_2 \geq 
\alpha$. Hence we may assume that  $b_1,b_2 <\alpha $. In this case $b_2+b_3+\cdots + b_d 
\geq \alpha +1$, and we can find $b'_2= b_2, b'_3 \leq b_3,...,b'_d \leq b_d$ such that 
$b'_2+\cdots + b'_d = \alpha $. Let $b'_1=0$. Then both elements $\bbf'$ and $\bbf-\bbf'$ 
satisfy the induction hypothesis, which implies $\bbf = \bbf'+ (\bbf-\bbf') \in S$.

Further, let $\bbf= (b_1,b_2,0,...,0)$ with $b_1+b_2 = \alpha (\alpha -2)$.  We also show that 
$\bbf\in S$. Indeed,  we can write $b_2 = p\alpha +q$, where $p\leq \alpha -2, \ q\leq  \alpha 
-1$. Note that $p= \alpha -2$ implies $q=0$ and $\bbf = (\alpha -2)\ebf_2 \in S$. Let $p\le 
\alpha -3$. In the case $\ p+q \geq  \alpha -1$, using also the equality
$ b_2  = (p+q - \alpha  +1 )\alpha  + (\alpha -q)(\alpha -1)$,
we may write $b_2 = m_1 + m_2(\alpha -1) + m_3\alpha $ with $m_1+m_2+m_3 \leq \alpha 
-2$. Hence 
$$\bbf = m_1\abf_1 + m_2 \abf_2 + m_3\ebf_2 + (\alpha -2-m_1-m_2-m_3)\ebf_1 \in S.$$
Summarizing the above arguments we get that $\bbf \in S$ if $\deg (\bbf) =\alpha -2$. 

Now let $\abf \in (\alpha -1)\Acal$. Since $\alpha \geq d+1$,   $a_1+ \cdots + a_d = \alpha 
(\alpha -1) \geq d\alpha $ and there is an index $i$ such that $a_i \geq \alpha $. Note that 
$\deg (\abf - \ebf_i) = \alpha -2$. By the above result $\abf - \ebf_i \in S$. Hence $\abf = \ebf_i 
+ (\abf-\ebf_i) \in \{\ebf_1,...,\ebf_d\} + S$, which implies $r(S) \leq \alpha -2$.

Summing up we get $r(S) = \alpha -2$. 

On the other hand,  let the term order be specified  as above.   Note that $x_1^{\alpha -1} - 
x_2y_1^{\alpha -2} \in I_\Acal$,  $x_1^{\alpha -1} > x_2y_1^{\alpha -2}$ and there is no 
other binomial of $I_\Acal$ whoes first term divides $x_1^{\alpha -1}$. Therefore the 
binomial $x_1^{\alpha -1} - x_2y_1^{\alpha -2}$ must be contained in the reduced 
Gr\"obner basis of $I_\Acal$. The degree of this binomial is $r(S)+1$, which is bigger than the 
half of the bound of Theorem \ref{A1}.
}\end{Example} 

In \cite{HS} some bounds for $r(S)$ were given. Let us recall them here. Let $\Pcal$ denote 
the convex polytope spanned by $\Acal \subset \Rset^d$. Note that $\Pcal$ is a
$(d-1)$-dimensional polytope whose faces are spanned by 
$$\Acal_I = \{\abf\in \Acal;\ a_{i} = 0\ \ \text{for\ all\ } i\in I\},$$
where $I \subseteq \{1,...,d\}$. Let
 $\Pcal_I$ denote the corresponding  face of $\Pcal$. (For short, we will also write $\Acal_i,\ 
\Pcal_i$ instead of $\Acal_{\{i\}},\ \Pcal_{\{i\}}$.)
  We  say that a face $\Pcal_I$ is {\it full} if $\Acal_I$ contains all points of 
$M_{\alpha, d}$ lying on this face, i.e. if $\Acal_I = \Pcal_I \cap M_{\alpha,d}$. By 
\cite[Theorem 1.1, Lemma 1.2 and Lemma 1.3]{HS} we have

\begin{Lemma} \label{A2} Let $\deg K[S]$ denote the multiplicity of $K[S]$. Then
\begin{itemize} \item[{\rm (i)}]  $r(S) \leq \deg K[S] - \codim K[S]$.
\item[{\rm (ii)}] If $\Pcal$ has a full face of dimension $i$, then
$r(S) \leq \alpha^{d-1-i} + i-1.$
\item[{\rm (iii)}] If  a $p$-dimensional face $\Pcal_I$ contains at least 
 $q+p+1$ points of $\Acal$, where $p\leq d-1$, then
 $r(S) \leq   (\alpha^p - q)\alpha^{d-1-p}.$
 \end{itemize}
 \end{Lemma}
 
 From Theorem \ref{A1} and Lemma \ref{A2}(i) we immediately get a bound which is very 
close to the Eisenbud-Goto bound.
 
 \begin{Theorem} \label{A3} With respect to any term order specified as above, $I_\Acal$ has 
a Gr\"obner basis of degree at most $\max\{ 2,\ 2(\deg K[S] - \codim K[S])-1)\}$.
 \end{Theorem}
 
 It should be noted that if $S$ is not necessarily a simplicial semigroup, then Sturmfels 
\cite{St2} showed that w.r.t. any term order, the ideal $I_\Acal$ has a Gr\"obner basis of 
degree at most $c\cdot \deg K[S]$.

The estimations in Theorem \ref{A1} and Theorem \ref{A3} may not hold if the term order 
does not satisfy the condition (*).
 
 \begin{Example} {\rm Let $\Acal =\{(4,0),(3,1),(1,3),(0,4)\}$. Then
 $$I_\Acal = ( x_1x_2-y_1y_2,x_1^3-x_2y_1^2,x_2^3 -x_1y_2^2, x_2^2y_1 - x_1^2 
y_2).$$
 The above minimal basis of $I_\Acal$ is also a minimal Gr\"obner basis w.r.t. the reverse 
lexicographic order. W.r.t. the lexicographic order we get the following minimal Gr\"obner 
basis:
 $$\{ x_1x_2-y_1y_2,x_1^3-x_2y_1^2, x_1y_2^2- x_2^3 ,  x_1^2 y_2 - x_2^2y_1, x_2^4 
-y_1y_2^3\}.$$
 In this example $r(S) = \deg K[S]- \codim K[S] =2$ and both bounds in Theorem \ref{A1} and 
Theorem \ref{A3} are equal to $3$.}
 \end{Example}
 
 Note that  a similar, but weaker bound was given for the Castelnuovo-Mumford regularity in 
\cite[Theorem 3.5]{HS}. There it was also shown  that the Castelnuovo-Mumford regularity of 
$K[S]$ is bounded by $c(\alpha -1)$ (see \cite[Theorem 3.2(i)]{HS}). It turns out that for 
certain Gr\"obner bases this bound also holds.
 
 \begin{Theorem} \label{A4} With respect to any term order specified as above, $I_\Acal$ has 
a Gr\"obner basis of degree at most $\max\{c, \alpha ,\ c(\alpha -1) -1\} \leq c(\alpha -1)$.
 \end{Theorem}
 
 \begin{pf} The proof is similar to that of Theorem \ref{A1}. Let $s= \max\{c, \alpha ,\ c(\alpha 
-1) -1\}$ and set
$$G =\{ \xbf^\mbf\ybf^\nbf - \xbf^\pbf\ybf^\qbf \in I_\Acal |\  \deg (\xbf^\mbf\ybf^\nbf ) = 
\deg( \xbf^\pbf\ybf^\qbf ) \leq s\}.$$
 Assume that $G$ is not a Gr\"obner basis. Then one can find a binomial $b= 
\xbf^\mbf\ybf^\nbf - \xbf^\pbf\ybf^\qbf\in I_\Acal$ of the smallest degree $\deg b>s$ such 
that $\gin(g) \nmid \xbf^\mbf\ybf^\nbf $ for all $g\in G$. Since $\alpha \abf_i = a_{i1}\ebf_1 
+ \cdots + a_{id}\ebf_d$, $x_i^\alpha - \ybf^{\abf_i} \in G$ for all $i=1,...,c$. Note that 
$x_i^\alpha > \ybf^{\abf_i}$. Since $\gin(x_i^\alpha - \ybf^{\abf_i}) \nmid 
\xbf^\mbf\ybf^\nbf $, We must have  $m_i \leq \alpha -1$ for all $i\leq c$.  

If $p_i \geq \alpha $, then
$$\xbf^\mbf\ybf^\nbf  - \frac{\xbf^\pbf}{x_i^\alpha}\ybf^{\qbf + \abf_i} = 
(\xbf^\mbf\ybf^\nbf - \xbf^\pbf\ybf^\qbf) + (x_i^\alpha - \ybf^{\abf_i}) 
\frac{\xbf^\pbf}{x_i^\alpha}\ybf^{\qbf} \in I_{\Acal}.$$
Note that $\xbf^\mbf\ybf^\nbf > \xbf^\pbf\ybf^\qbf > \frac{\xbf^\pbf}{x_i^\alpha}\ybf^{\qbf 
+ \abf_i}$. Replacing $b$ by $\xbf^\mbf\ybf^\nbf  - \frac{\xbf^\pbf}{x_i^\alpha}\ybf^{\qbf + 
\abf_i}$ and repeating this procedure, we may also assume that $p_i \leq \alpha -1$ for all 
$i\leq c$.  

As in the proof of Theorem \ref{A1}, let $C= \{ i|\ m_i \neq 0\}$ and $D= \{j|\ n_j \neq 0\}$. 
Then  we can also conclude that 
\begin{equation} \label{EA4}
\sum_{j\in D} n_j \leq \sum_{i\not\in C} p_i \leq (c - \sharp C) (\alpha -1),
\end{equation}
and that $\sum_{j\in D} n_j  = (c - \sharp C) (\alpha -1)$ implies  $\xbf^\mbf - \ybf^\qbf \in 
I_\Acal$.  Hence
$$\deg (\xbf^\mbf\ybf^\nbf ) = \sum_{i\in C}m_i + \sum_{j\in D} n_j \leq  \sharp C (\alpha 
-1) + (c - \sharp C) (\alpha -1) = c(\alpha -1).$$
Since $\deg (\xbf^\mbf\ybf^\nbf ) = \deg(b) \geq c(\alpha -1)$, we must have $\deg 
(\xbf^\mbf\ybf^\nbf ) = c(\alpha -1)$. Therefore $\sum_{j\in D} n_j  = (c - \sharp C) (\alpha 
-1)$ and $m_i  = \alpha -1$ for all $i\in C$.  By (\ref{EA4}) we have  $\xbf^\mbf - \ybf^\qbf 
\in I_\Acal$. If $C \neq \{1,...,c\}$, then $\deg(\xbf^\mbf ) \leq s$ and $\xbf^\mbf - \ybf^\qbf 
\in G$, which is impossible because $\xbf^\mbf \mid \gin(b)$. Thus $C= \{1,...,c\}$. This 
yields $D= \emptyset $ and
$$b = (x_1\cdots x_c)^{\alpha -1} - \ybf^\qbf.$$
Let $\abf = \abf_1 + \cdots + \abf_c$. The above equality assures that $\alpha \mid (\alpha 
-1)a_i $ for all $i=1,...,c$. This implies $a_i = q'_i\alpha $ for some $q'_i\in \Nset$. But then 
$g:= x_1\cdots x_c - y_1^{q'_1}\cdots y_d^{q'_d} \in I_\Acal$. Since $\deg(x_1\cdots x_c) 
= c \leq s$, $g\in G$ and we get a contradiction that $\gin(g) =  x_1\cdots x_c\mid \gin(b) = ( 
x_1\cdots x_c)^{\alpha -1}$. The proof of the theorem is completed. \hfill $\square$
\end{pf}

The following example shows that the bound in the above theorem is not very far from the best 
what we can hope.

\begin{Example} \label{A5} {\rm Let $d=c \geq 3, \ \alpha = 2\beta $, where $\beta \geq 2$, 
and
$$\begin{array}{ll}
\abf_1 & = (1, 2\beta -1, 0,...,0),\
\abf_2 = (0, 1, 2\beta -1, 0,...,0),\ ...\\
\abf_{c-1} & = (0,...,0,1,2\beta -1),\
\text{and}\ \abf_c = (\beta ,0,...,0,\beta ).
\end{array}$$
Then
$$g:= x_1^\beta \cdots x_{c-1}^\beta - x_cy_2^\beta \cdots y_{c-1}^\beta y_c^{\beta -1} \in 
I_\Acal.$$
Assume that there is a binomial
$$b= x_1^{n_1}\cdots x_{c-1}^{n_{c-1}} - x_1^{p_1}\cdots x_c^{p_c} y_1^{q_1}\cdots 
y_d^{q_d} \in I_\Acal,$$
where $n_1,...,n_{c-1} \leq \beta ,\ p_1 <2\beta ,..., p_c< 2\beta $ and if $n_i \neq 0$ then 
$p_i=0$. Then
\begin{equation} \label{EA5}
n_1\abf_1 + \cdots + n_{c-1}\abf_{c-1} = p_1\abf_1 + \cdots + p_c\abf_c + q_1\ebf_1 + 
\cdots + q_d\ebf_d.
\end{equation}
Comparing the first coordinates of both sides, we obtain 
\begin{eqnarray}
n_1 =& p_1 + \beta p_c + 2\beta q_1, \label{EA5a}\\
-n_1+n_2 =& -p_1 + p_2 + 2\beta (p_1+q_2 - n_1) ,\label{EA5b}\\
&... \nonumber \\
-n_{c-2} + n_{c-1} =& -p_{c-2} + p_{c-1} + 2\beta (p_{c-2} + q_{c-1} - n_{c-2}).\nonumber 
\end{eqnarray}
Using the property that  $p_i=0$ if $n_i \neq 0$, the above relations imply that $\beta \mid 
n_i, p_i $ for all $i<c$. Hence for each $i$, either $n_i=\beta $ or $n_i=0$ and also either 
$p_i=\beta $ or $p_i=0$. If $n_1 =0$, by (\ref{EA5a}) it would give $p_1= p_c=0$. Let $i\leq 
c-1$ be the least index such that $n_i\neq 0$. Then $p_1=\cdots = p_i =0$. A comparision of 
the $i$-th coordinates would give $n_i = 2\beta q_i$, a contradiction. Thus $n_1 = \beta $.  
This implies by (\ref{EA5a}) that $p_1=0$ and $p_c=1$. If  $n_2=0$, then by (\ref{EA5b}) 
we would get $p_2 = \beta $. The third coordinate in the left hand  of (\ref{EA5}) would be at 
most $\beta $ (a possible contribution of $n_3\abf_3$ when $c\geq 4$), while that in the right 
hand would be at least  $\beta (2\beta -1) + 2\beta q_3$, a contradiction. Therefore 
$n_2=\beta , p_2=0$. Repeating this procedure, at the end we get  $n_1=\cdots = n_{c-1} = 
\beta $, and so $b= g$. This means  the binomial $g$ defined above should be an element of 
the reduced Gr\"obner basis of $I_\Acal$ w.r.t. any term order specified as above. We have 
$\deg (g) = (c-1)\beta $, while the bound of Theorem \ref{A4} is $c(2\beta -1)$.
}\end{Example}

 The last result of this section shows that, by using a suitable term order, the computation of 
simplicial toric ideals runs rather quickly.
  In order to compute $I_\Acal$, a standard procedure is the following (see, e.g., \cite{St1}, 
Algorithm 4.5):
  
  \begin{itemize}
  \item[1.] Form the ideal $J_\Acal = (x_1-\tbf^{\abf_1},...,x_c-\tbf^{\abf_c}, y_1-t_1^\alpha 
,...,y_d-t_d^\alpha ) \subset K[\tbf,\xbf,\ybf]$.
  \item[2.] Compute a Gr\"obner basis $G'$ of $J_\Acal$ by Buchberger's algorithm, using an 
elimination order $\preceq $ with respect to the variables $t_1,...,t_d$.
  \item[3.] From $G'$ get a Gr\"obner basis $G= G' \cap K[\xbf,\ybf]$ of $I_\Acal = J_\Acal 
\cap K[\xbf,\ybf]$.
  \end{itemize}
  
  \begin{Proposition}\label{A6} Assume that the restriction of the elimination order $\preceq 
$ on $K[\xbf,\ybf]$ either is the reverse lexicographic order or satisfies the condition (*). Then 
$J_\Acal$ has a Gr\"obner basis of degree at most $d(\alpha -1) + \min\{2r(S) , c(\alpha -1)\}$.
  \end{Proposition}
  
  \begin{pf} The proof is similar to that of Theorem \ref{A1}. Let $s= d(\alpha -1) + \min\{2r(S) 
, c(\alpha -1)\}$ and set
$$G =\{ \tbf^\pbf\xbf^\mbf\ybf^\nbf - \tbf^{\pbf'}\xbf^{\mbf'}\ybf^{\nbf'} \in J_\Acal |\  \deg 
(\tbf^\pbf\xbf^\mbf\ybf^\nbf ) = \deg( \tbf^{\pbf'}\xbf^{\mbf'}\ybf^{\nbf'} ) \leq s\}.$$
 Assume that $G$ is not a Gr\"obner basis. Then one can find a binomial $b= 
\tbf^\pbf\xbf^\mbf\ybf^\nbf - \tbf^{\pbf'}\xbf^{\mbf'}\ybf^{\nbf'}\in J_\Acal$ of the smallest 
degree $\deg b>s$ such that 
$\gin(g) \nmid \tbf^\pbf \xbf^\mbf\ybf^\nbf $ for all $g\in G$. Since $t_i^\alpha - y_i\in G$ 
and $t_i^\alpha \succ  y_i$,  we must have $p_i\leq \alpha -1$ for all $i\leq d$. Let
$$G_1 = \{ \xbf^\mbf\ybf^\nbf - \xbf^\pbf\ybf^\qbf \in I_\Acal |\  \deg (\xbf^\mbf\ybf^\nbf ) = 
\deg( \xbf^\pbf\ybf^\qbf ) \leq 2r(S)\},$$
and
$$G_2 = \{ \xbf^\mbf\ybf^\nbf - \xbf^\pbf\ybf^\qbf \in I_\Acal |\  \deg (\xbf^\mbf\ybf^\nbf ) = 
\deg( \xbf^\pbf\ybf^\qbf ) \leq c(\alpha -1)\}.$$
Then $ G_1\cup G_2\subset  G$. Assume that $\deg(\xbf^\mbf\ybf^\nbf) > 2r(S)$. Under
 the assumption on the term order,  it was shown in the proof of   Theorem \ref{A1} that 
$\deg(\xbf^\mbf\ybf^\nbf)\ \vdots\ \gin(g)$ for some $g\in G_1$. Then $\gin(g) \mid \tbf^\pbf 
\xbf^\mbf\ybf^\nbf $ which is impossible. Hence $\deg(\xbf^\mbf\ybf^\nbf) \leq 2r(S)$. 
Similarly, using $G_2$ and analyzing the proof of  Theorem \ref{A4} we also get 
$\deg(\xbf^\mbf\ybf^\nbf) \leq c(\alpha -1)$. Therefore, 
$$\deg (b)= \sum p_i + \deg(\xbf^\mbf\ybf^\nbf) \leq d(\alpha -1) + \min\{2r(S) , c(\alpha 
-1)\}= s,$$ a contradiction. \hfill $\square$
 \end{pf}

\section{Eisenbud-Goto bound}\label{B}

In this section we will provide some partial positive answers to the question posed in the 
introduction. As a fisrt result we have the following immediate consequence of Theorem 
\ref{A4} (cf. \cite[Corollary 3.6]{HS}  for the Castelnuovo-Mumford regularity):

\begin{Corollary} \label{B1} Assume that $\codim K[S] \leq \deg K[S]/\alpha $. Then, w.r.t. 
any term order specified in Section \ref{A}, $I_\Acal$ has a Gr\"obner basis of degree at most 
$\deg K[S] - \codim K[S] $.
\end{Corollary}

\begin{Remark} {\rm By \cite{HS}, $\alpha \mid \deg K[S]$. Very often we have $\deg K[S] = 
\alpha^{d-1}$ (see \cite[Example 1.4]{HS}). Thus, if $d\geq 3$ and the codimension is not too 
big, the Eisenbud-Goto bound holds for a Gr\"obner basis of $I_\Acal$.
Unfortunately it may happen that $\alpha = \deg K[S]$. Hence the above corollary even does 
not completely solve the case $c=2$. This case was solved by Peeva and Sturmfels  by 
considering a more general class - the class of codimension 2 lattice ideals (see \cite[Theorem 
7.3 and Proposition 8.3]{PS}). Another proof was recently given in \cite{BGM} (see Theorems 
2.1, 2.8 and 3.5 there).
}\end{Remark}

Recall that a quotient ring $R/I$ modulo an homogeneous  ideal $I$ is said to be a {\it 
generalized Cohen-Macaulay} ring if all local cohomology modules $H^i_\mfr(R/I),\ i<\dim 
R/I,$ with the support in the maximal homogeneous ideal $\mfr$ of $R/I$ are of finite length 
(see the Appendix in \cite{SV}). The Castelnuovo-Mumford regularity of a finitely generated 
graded $R$-module $M$ is the number 
$$\reg(M) = \max\{ n|\ [H^i_\mfr(M)]_{n-i} \neq 0 \ \text{for} \ i\geq 0\}.$$
Note that $\reg(I) = \reg(R/I) + 1$. The following result is a simple observation, but has some 
interesting consequences.

\begin{Lemma} \label{C1} Assume that $K[S]$ is a generalized Cohen-Macaulay ring. Then 
$I_\Acal$ has a Gr\"obner basis w.r.t. the reverse lexicographic order of degree at most $\reg 
I_\Acal$.
\end{Lemma}

\begin{pf} Note that $y_1,...,y_d$ is a system of parameters of $K[S]$. Since $K[S] \cong  
K[\xbf]/I_\Acal$ is a generalized Cohen-Macaulay ring, the ideal $I_\Acal$ and all ideals 
$(I_\Acal, y_d,...,y_i),\ i=d,d-1,...,1$, are unmixed  up to  $\mfr$-primary components (see 
\cite[Proposition 3 in the Appendix]{SV}). In particular, $y_{i-1}$ is a non-zero divisor on the 
ring $K[\xbf]/(I_\Acal, y_d,...,y_i)^{\rm sat}$, where $J^{\rm sat}= \cup_{n\geq 1} J:\mfr^n$ 
denotes the saturation of $J$. This means $y_d,...,y_1$ is a generic sequence of $K[S]$ in the 
sense of \cite[Definition 1.5]{BS}. By \cite[Corollary 2.5]{BS},  the maximal degree of a 
minimal Gr\"obner basis of $I_\Acal$ is bounded by $\reg(I_\Acal)$. \hfill $\square$
\end{pf}

\begin{Remark}\label{C1b} {\rm  (i) The generalized Cohen-Macaulay property of $k[S]$ can 
be characterized combinatorically. For this purpose, let $S_i = \{\abf \in S|\ a_i =0\}$, and $S' = 
\cap_{i=1}^c (S-S_i)\supseteq  S$. Then $K[S]$ is a generalized Cohen-Macaulay ring if and 
only if the set $S'\setminus S$ is finite. The necessary part follows from \cite[Corollary 
2.3]{TH}. The sufficient part follows from the fact that $K[S']$ is always a Cohen-Macaulay 
ring (see \cite[Theorem 5.1]{GSW}  or \cite[Corollary 4.4]{TH}) and the following exact 
sequence:
$$0 \rightarrow K[S] \rightarrow K[S'] \rightarrow K[S'\setminus S] \rightarrow 0.$$

(ii) With respect to any term order specified in Section \ref{A}, $y_d,...,y_1$ is a system of 
parameters of $\gin (I_\Acal)$. This follows from the fact that $x_i^\alpha \in \gin (I_\Acal)$ 
for all $i\leq c$ (since $x_i^\alpha - \tbf^{\abf_i} \in I_\Acal$). However, if   $K[S]$ is not a 
generalized Cohen-Macaulay ring, this sequence maybe no more  a generic sequence of 
$K[S]$. For example, consider a simple case: $d=\alpha =3$ and
$$\begin{array}{ll}  \Acal = \{&\ebf_1,\ebf_2,\ebf_3, \abf_1= (2,0,1),\ \abf_2 = (1,2,0), 
\abf_3 = (1,1,1),\\ & \abf_4=(1,0,2), \abf_5=(0,2,1),\ \abf_6=(0,1,2)\}.\end{array}$$
Then w.r.t. the reverse lexicographic order 
$$\begin{array}{ll} \gin(I_\Acal) = (&x_1x_2, x_2x_3, x_2x_5, x_1^2, x_1x_3, x_3^2, 
x_2x_4, x_2x_6, x_3x_5, x_5^2, x_1x_4, x_3x_4,\\ & x_4x_5, x_4^2, x_3x_6, x_5x_6, 
x_4x_6, x_6^2, x_2^3, x_1x_6y_2) .
\end{array}$$
 Clearly $(\gin(I_\Acal),y_3)^{sat} = (\gin(I_\Acal),y_3)$ and $y_2$ is a zero divisor of 
$K[\xbf]/(\gin(I_\Acal), y_3)$. Hence, by \cite[Theorem 2.4(a)]{BS}, $y_3,y_2,y_1$ is not a 
generic sequence of  $K[\xbf]/I_\Acal$.
Note that in this example $S'\setminus  S= \langle M_{3,3}\rangle \setminus  S = (2,1,0) + 
\Nset\ebf_1$ is infinite.}
\end{Remark}
By \cite{GLP} we know that $\reg(I_\Acal) \leq \deg K[S] - \codim K[S]+1$ if $d=2$. Hence, 
by Lemma \ref{C1} we get

\begin{Corollary} \label{C3} Assume that $d=2$. Then, w.r.t. the reverse lexicographic order, 
$I_\Acal$ has a Gr\"obner basis of degree at most $\deg K[S] - \codim K[S]+1$.
\end{Corollary}

\begin{Corollary} \label{C4} Let $K[S]$ be a simplicial semigroup ring with isolated 
singularity. Then, w.r.t. the reverse lexicographic order,
\begin{itemize}
\item[{\rm (i)}] $I_\Acal$ has a Gr\"obner basis of degree at most $(d-1)(\alpha -2) +1$;
\item[{\rm (ii)}] $I_\Acal$ has a Gr\"obner basis of degree at most $\deg K[S] - \codim 
K[S]+1$.
\end{itemize}
\end{Corollary}

\begin{pf} Under the assumption, $K[S]$ is a generalized Cohen-Macaulay ring. On the other 
hand, by \cite[Theorem 2.1 and Corollary 2.2]{HH},  we have $\reg (K[S] )$ $ \leq (d-1)(\alpha 
-2)$ and $\reg (K[S] ) \leq \deg K[S] - \codim K[S]$. Since $\reg (I_\Acal) = \reg (K[S]) + 1$, 
the statement follows from Lemma \ref{C1}. \hfill $\square$
\end{pf}

\begin{Remark} \label{C5} {\rm By (4) in the proof of \cite[Theorem 2.1]{HH},  the condition 
of Corollary \ref{C4} is equivalent to the condition that $\Acal$ contains all points of 
$M_{\alpha ,d}$ of type $(0,..,\alpha -1,...,1,...,0)$, where $\alpha -1,1$ stay in the $i$-th and 
$j$-th positions, respectively, and the other coordinates are zero.}
\end{Remark}

Recall that a (not necessarily simplicial) semigroup $S$ is said to be {\it normal} if $S = \{ \abf 
\in \Zset(S) |\ n\abf \in S \ \text{for some}\  n>0\}$, where $\Zset(S)$ denotes the subgroup of 
$\Zset^d$  generated by $S$. Under this condition, it is well-known  that $\reg(I_\Acal) \leq 
d$ (this holds even without the assumption $S$ being simplicial, see \cite[Proposition 
13.14]{St1}). From this, Sturmfels posed the following question:

 {\bf Question} (\cite[p. 136]{St1}):  {\it If the semigroup $S$ is normal, does the toric ideal 
$I_\Acal$ posses a Gr\"obner basis of degree at most $d$}?

In the above problem, the semigroup $S$ is not necessarily simplicial. It turns out that the 
simplicial case immediately follows from Lemma \ref{C1}.

\begin{Corollary} \label{C2} Assume that the simplicial semigroup $S$ is normal. Then the 
toric ideal $I_\Acal$ has a Gr\"obner basis w.r.t. the reverse lexicographic order of degree at 
most $d$.
\end{Corollary}
\begin{pf} When $S$ is normal, $K[S]$ is a Cohen-Macaulay ring by \cite{Ho}. Hence  the 
statement follows from Lemma \ref{C1} and the above mentioned result on $\reg(I_\Acal)$. 
\hfill $\square$
\end{pf}

\begin{Remark} \label{C2b} {\rm Let $S$ be a simplicial normal semigroup. In  general $d$ 
is much smaller than the Eisenbud-Goto bound. However, even if this is not the case, 
$I_\Acal$ still  has a Gr\"obner basis w.r.t. the reverse lexicographic order of degree at most 
$\deg K[S] - \codim K[S]+1$. This follows from Lemma \ref{C1} and the fact that the 
Eisenbud-Goto bound holds  for the regularity of  perfect prime ideals (see \cite{Tr}).}
\end{Remark}

Under the assumption of the following result it was shown in \cite[Proposition 3.7]{HS} that 
$\reg(I_\Acal) \leq \deg K[S] - \codim K[S]+1$. Unfortunately we cannot use it to derive the 
corresponding result for a Gr\"obner basis, because $K[S]$ is in general not a generalized 
Cohen-Macaulay ring.

\begin{Proposition} \label{B2} Assume that $\deg K[S] = \alpha^{d-1}$ and $\alpha \leq 
d-1$. Then, w.r.t. the reverse lexicographic order, $I_\Acal$ has a Gr\"obner basis of degree at 
most $\deg K[S] - \codim K[S]+1$.
\end{Proposition}

\begin{pf} If $\Acal = M_{\alpha ,d}$, then $K[S]$ is normal and by Remark \ref{C2b} we are 
done.  Hence we may assume that
$$c\leq \sharp M_{\alpha,d}-1-d = {\alpha +d-1\choose d-1}-d-1.$$
If $\alpha \geq 3,\ d\geq 6$ or $\alpha =4,\ d=5$, then by \cite[Claim 1]{HS},  $c\leq 
\alpha^{d-2}$. Hence
$$\deg K[S] - c+1 > \alpha^{d-1} - \alpha^{d-2} = \alpha^{d-2}(\alpha -1) \geq c(\alpha -1),$$
and by Theorem \ref{A4} we are done. Thus the left cases are: $\alpha =2, \ d\geq 3$; $\alpha 
=3,\ d=4$ and $\alpha =3,\ d=5$. We consider these cases separately.
\vskip0.3cm

{\bf Case 1}: $\alpha =2, \ d\geq 3$. Then $c\leq d(d+1)/2 - (d+1) =(d-2)(d+1)/2$. It is easy to 
verify that $(d-2)(d+1)/2-1\leq 2^{d-2}$. Hence, if $c\leq (d-2)(d+1)/2-1$ we have $c\leq 
2^{d-1} - c = \deg K[S] -c$. By Theorem \ref{A4} we are done. The left case is 
$c=(d-2)(d+1)/2$, i.e. $\Acal$ is obtained from $M_{2,d}$ by deleting exactly one point. We 
may assume $\Acal = M_{2,d}\setminus \{\bbf = (1,1,0,...,0)\}$. Note that $2\abf_i \in 
\{\ebf_1,...,\ebf_d\} + \{\ebf_1,...,\ebf_d\} $ for all $i\leq c$. Moreover, if $\abf_i,\abf_j$ are 
two differents points and  $\abf_i,\abf_j,\bbf$ do not lie in the same 2-dimensional face of 
$\Pcal$, then $\abf_i + \abf_j \in \Acal + \{\ebf_1,...,\ebf_d\} $ (see Fig. 1). From this it 
follows that $r(S)=2$. By Theorem \ref{A1}, $I_\Acal$ has a 
Gr\"obner basis of degree at most $3 \leq 2^{d-1} - (d-2)(d+1)/2 + 1 = \deg K[S] -c +1$. 

\centerline{\setlength{\unitlength}{0.4cm}
 \begin{picture}(8,7)
\put(0,2){\line(3,4){3}}
\put(-0.8,2){$\ebf_2$}
\put(1.5,4){\circle{0.2}} \put(1,4){$\bbf$}
\put(0,2){\line(3,-2){3}}
\put(3.1,6){$\ebf_1$}
\multiput(0,2)(0.3,-0.05){20}{.}
\put(3,0){\line(0,1){6}}
\put(2.8,-0.5){$\ebf_3$}
\put(3,0){\line(3,1){3}}
\put(3,6){\line(3,-5){3}}
\put(6.1,1){$\ebf_4$}
\put(3,3){\circle*{0.2}} \put(3.2,3){$\abf_k$}
\put(4.5,0.5){\circle*{0.2}} \put(4.5,-0.1){$\abf_i$}
\put(4.5,3.5){\circle*{0.2}} \put(4.7,3.5){$\abf_j$}
\put(1,-2){ Fig. 1}
\end{picture} 
\begin{picture}(11,7)
\put(3,0){\line(1,0){6}}
\put(3,0){\line(1,2){3}}
\put(6,6){\line(1,-2){3}}
\put(6,6){$\ebf_1$} \put(9,0){$\ebf_3$} \put(2.2,0){$\ebf_2$}
\put(5,0){\circle{0.2}} 
\put(7,0){\circle*{0.2}} \put(6.8,-0.6){$\abf_3$}
\put(4,2){\circle{0.2}} 
\put(5,4){\circle*{0.2}} \put(4,4){$\abf_1$}
\put(8,2){\circle*{0.2}} \put(8,2){$\abf_2$}
\put(7,4){\circle{0.2}} 
\put(6,2){\circle*{0.2}} \put(6.2,2){$\bbf_4$}
\put(4,-2){ Fig. 2}
\end{picture}}
\vskip1cm  

{\bf Case 2}: $\alpha = 3, d=4$. Then $\deg K[S] = 27$ and $c\leq 15$. By Theorem \ref{A4}, 
the statement of the proposition holds true for $c\leq 9$. Let $c\geq 10$, i.e. $\Acal$ is 
obtained from $M_{3,4}$ by deleting at most 6 points. We distinguish two subcases.
\vskip0.3cm

{\it Subcase 2a}: Each edge of $\Pcal$ contains exactly one deleting point. In this case 
$c=10$. By Theorem \ref{A4} it suffices to show that $r(S) \leq 8$.

Consider, for example, the facet $\Pcal_4 = \{\abf \in \Pcal|\ a_4 =0\}$. Then $\Acal_4 = \Acal 
\cap \Pcal_4$ has exactly 7 points, say
$\Acal_4 = \{\ebf_1, \ebf_2,\ebf_3, \abf_1,\abf_2,\abf_3,\bbf_4 = (1,1,1,0)\}$ as shown in 
Fig. 2, where $\abf_1$ can be taken as $(2,1,0,0)$, while there are two choices for each of 
$\abf_2$ and $\abf_3$. One can check by computer that in this case the reduction number 
$r(\langle \Acal_4\rangle) \leq 3$. In particular,  $\sum_{i=1}^3m_i\abf_i  + 
n_4\bbf_4\not\in  S + \{\ebf_1,\ebf_2+\ebf_3\}$ implies that $m_i,n_4\leq 2$ and
\begin{equation}\label{EB2}
\sum_{i=1}^3m_i + n_4 \leq 3.
\end{equation}
 Moreover, since
$2\bbf_4 + \abf_1 = (4,3,2,0)  = \ebf_2 + 2(2,0,1,0) = \ebf_1 + \ebf_2 + (1,0,2,0)$ and one of 
two points $(2,0,1,0) $ and $(1,0,2,0)$ on the edge $\overline{\ebf_1\ebf_3}$ must belong to 
$\Acal_4$, we get that $2\bbf_4 + \abf_1 \in S + \{\ebf_1,\ebf_2+\ebf_3\}$. The same is true 
for $2\bbf_4 +\abf_2$ and $2\bbf_4 +\abf_4$. This means, in addition to (\ref{EB2}) we also 
have $m_1=m_2=m_3 =0$ if $n_4=2$.

Finally, we can write $\Acal = \{\ebf_1, \ebf_2,\ebf_3, \abf_1,...,\abf_6,\bbf_1,...,\bbf_4 \},$ 
where $\bbf_i$ is the inner point of the facet $\Pcal_i$. Assume that 
$$\sum_{i=1}^6m_i\abf_i  + \sum_{j=1}^4 n_j\bbf_j \not\in  S + 
\{\ebf_1,\ebf_2,\ebf_3,\ebf_4 \}.$$
Then inequalities of Type (\ref{EB2}) should hold for all facets of $\Pcal$. Adding all these 
inequalities we get
$$\sum_{i=1}^6m_i + \frac{1}{2} \sum_{j=1}^4 n_j \leq 6.$$
If  $n_j \leq 1$ for  $j$, then $\sum_{i=1}^6m_i +  \sum_{j=1}^4 n_j \leq 8$. If, say $n_4$, is 
equal to $2$, then $m_1=m_2=m_3 =0$. Adding the inequalities on the facets $\Pcal_1, \ 
\Pcal_2$ and $\Pcal_3$ we will get $\sum_{i=4}^6m_i + \frac{1}{2} \sum_{j=1}^3 n_j \leq 
4\frac{1}{2}.$ Continuing as above, we also get in all cases that $\sum_{i=1}^6m_i +  
\sum_{j=1}^4 n_j \leq 8$. This means $r(S) \leq 8$.
\vskip0.3cm

{\it Subcase 3b}: At least one edge of $\Pcal$ is full. By Lemma \ref{A2}(ii), $r(S)\leq 9$. 
Hence, by Theorem \ref{A1}, the statement holds true if $c\leq 11$.  Moreover, if $\Pcal$ has 
a full facet, then again by Lemma \ref{A2}(ii), $r(S)\leq 5$, and by Theorem \ref{A1} we are 
done. Hence, we may assume that $c=12,13,14$, and $\Pcal$ has no full facet. This 
corresponds to the situation when $\Pcal$ has 2,3 or 4 deleting points.

Assume that $\Pcal$ has a facet, say $\Pcal_4$, which contains exactly one deleting point 
$\bbf$, i.e. one can write $\Acal_4 = \{\ebf_1, \ebf_2,\ebf_3,\abf_1,...,\abf_6\}$ and 
$\bbf\not\in \Acal_4$. By Theorem \ref{A1}, it suffices to show that $r(S)\leq 5$. If this is not 
the case, then one can find $m_1,...,m_c\in \Nset$ such that $\sum_{i=1}^c m_i=6$ and 
\begin{equation} \label{EB2b}
\sum_{i=1}^c m_i\abf_i \not\in S+ \{\ebf_1,\ebf_2,\ebf_3,\ebf_4\}.
\end{equation}
We follow the idea in the proof of \cite[Lemma 1.2]{HS}. Considering $6$ subsums 
$\abf_1,...,m_1\abf_1, \ m_1\abf_1+\abf_2,...,\ \sum_{i=1}^c m_i\abf_i$ we can find either 
two subsums whose last coordinates are divisible by $3$, or three subsums whose last 
coordinates are congruent modulo $3$. Taking also the differences of these subsums, we can 
find in both cases two subsums $\bbf_1 = \sum p_i\abf_i$ and $\bbf_2 = \sum q_i \abf_i$ 
such that $p_i+q_i \leq m_i,\ i\leq 3$, $\deg(\bbf_1),\deg(\bbf_2) \geq 2$ and the last 
coordinates of $\bbf_1,\bbf_2$ are divisible by $3$.  Hence we can write $\bbf_i = \bbf_i' + 
n_i \ebf_4,\ i=1,2$, where $\bbf_i' \in \langle \Acal_4, \bbf\rangle$ . By (\ref{EB2b}) we 
must have $\bbf_1,\bbf_2 \not\in S+ \{\ebf_1,\ebf_2,\ebf_3,\ebf_4\}$, which yields $0 \neq 
\bbf_1',\bbf_2' \not\in \{\ebf_1,\ebf_2,\ebf_3\} + \langle \Acal_4\rangle$. Together with the 
fact  $2\bbf \in \langle \Acal_4\rangle$, this implies $\bbf_i' \in \{\bbf, \abf_1,...,\abf_6\} + 
\langle \Acal_4\rangle$. Since also all elements $2\bbf, \bbf+\abf_1,...,\bbf+\abf_6 \in \langle 
\Acal_4\rangle$, the previous relation assures that $\bbf_1'+\bbf_2'\in \langle \Acal_4\rangle 
\subset S$. By (\ref{EB2b}) we must have $n_1=n_2 =0$, and so $\bbf_1+\bbf_2\in \langle 
\Acal_4\rangle$. However it is easy (or using computer) to see  that $r(\langle \Acal_4\rangle) 
=2$. Since $\deg(\bbf_1+\bbf_2)\geq 4$, 
$$\bbf_1+\bbf_2\in \langle \Acal_4\rangle + \{\ebf_1,\ebf_2,\ebf_3\} \subseteq S+ 
\{\ebf_1,\ebf_2,\ebf_3,\ebf_4\},$$
which contradicts (\ref{EB2b}) because $p_i+q_i \leq m_i,$ for all $ i\leq 3$.

Thus, each facet of $\Pcal$ must have at least two deleting points. In particular, $c=12$ and 
$\Pcal$ has exactly 4 deleting points. There are only two situations shown in Fig. 3 and Fig. 4. 
In the situation of Fig. 3 there is eventually one configuration, and by computer we see that 
$r(S)=2$. In the situation of Fig. 4 one can show as in Subcase 2a (or using computer for eight 
different configurations), that $r(S)\leq 8$. But then by Theorem \ref{A1}, $I_\Acal$ has a 
Gr\"obner basis of degree at most $15 < \deg K[S] - c+ 1= 16$. The Subcase 2b is completely 
solved.
 
\centerline{\setlength{\unitlength}{0.4cm}
 \begin{picture}(8,7)
\put(0,2){\line(3,4){3}}
\put(1,3.3){\circle{0.2}} 
\put(2,4.7){\circle{0.2}}
\put(0,2){\line(3,-2){3}}
\multiput(0,2)(0.3,-0.05){20}{.}
\put(3,0){\line(0,1){6}}
\put(3,0){\line(3,1){3}}
\put(3,6){\line(3,-5){3}}
\put(4,0.3){\circle{0.2}} 
\put(5,0.7){\circle{0.2}} 
\put(1,-2){ Fig. 3}
\end{picture} 
\begin{picture}(8,7)
\put(0,2){\line(3,4){3}}
\put(1,3.3){\circle{0.2}} 
\put(3,4){\circle{0.2}}
\put(0,2){\line(3,-2){3}}
\multiput(0,2)(0.3,-0.05){20}{.}
\put(3,0){\line(0,1){6}}
\put(3,0){\line(3,1){3}}
\put(3,6){\line(3,-5){3}}
\put(4,0.3){\circle{0.2}} 
\put(2,1.7){\circle{0.2}} 
\put(1,-2){ Fig. 4}
\end{picture} }
\vskip1cm

{\bf Case 3}: $\alpha =3, d=5$. We have $c\leq 29$ and $\deg K[S] = 81$. If $c\leq 27$, then 
$\deg K[S] - c + 1 \geq 54 \geq 2c$, and by Theorem \ref{A4} we are done. If $c= 28, 29$, 
then $\Acal$ is obtained from $M_{3,5}$ by deleting 1 or 2 points. But then $\Pcal$ has a 
full $2$-dimensional face. By Lemma \ref{A2}(ii), $r(S) \leq 10$. Hence, by Theorem 
\ref{A1}, we are  also done in this subcase. \hfill $\square$
\end{pf}

Finally we show that if on an edge of $\Pcal$ there are enough points belonging to $\Acal$, 
then the Eisenbud-Goto bound also holds for the maximal degree of a Gr\"obner basis of 
$I_\Acal$. Note that in this setting, the Eisenbud-Goto conjecture on $\reg(I_\Acal)$  is still 
not verified (cf. \cite[Corollary 3.8]{HS}).

\begin{Proposition} \label{B3} Assume that $\deg K[S]= \alpha^{d-1}$ and there exists an 
edge of $\Pcal$ such that it is either full or at least $(\frac{3}{4} + \frac{1}{4d})\alpha +2$ 
integer points on it belong to $\Acal$. Then w.r.t. the reverse lexicographic order, $I_\Acal$ 
has a Gr\"obner basis of degree at most $\deg K[S] - \codim K[S] +1$.
\end{Proposition}

\begin{pf} By Corollary \ref{C3}, Remark \ref{C2b} and Proposition \ref{B2} we may assume 
that $\alpha \geq d\geq 3$ and
$$c\leq {\alpha + d-1 \choose d-1} - d-1.$$
First we consider the case when at least $(\frac{3}{4} + \frac{1}{4d})\alpha +2$ integer 
points on an edge belong to $\Acal$. By Lemma \ref{A2}(iii), $r(S) \leq 
\frac{d-1}{4d}\alpha^{d-1}$. Hence, by Theorem \ref{A1}, it suffices to show that
$$\alpha^{d-1} - {\alpha + d-1 \choose d-1} + d+1 \geq \frac{d-1}{2d}\alpha^{d-1} -1,$$
or equivalently
\begin{equation}\label{EB3}
\frac{d+1}{2d}\alpha^{d-1} + d+2 \geq {\alpha + d-1 \choose d-1}.
\end{equation}
We show this by induction on $d\geq 3$. For $d=3$ this is equivalent to $\alpha^2 - 9\alpha 
+24\geq 0$. So, assume that the inequality holds for $d\geq 3$. In the dimension $d+1$, by 
induction we have
$$\begin{array}{ll} 
 \displaystyle{{\alpha + d \choose d}} &=  \displaystyle{ \frac{\alpha +d}{d}{\alpha + d-1 
\choose d-1}
 \leq \frac{\alpha +d}{d}(\frac{d+1}{2d}\alpha^{d-1} + d+2) }  \\
&=  \displaystyle{ \frac{(\alpha +d)(d+1)}{2d^2}\alpha^{d-1} + \frac{d+2}{d}\alpha + d+2.}
\end{array}$$
Hence
$$\begin{array}{l} 
\displaystyle{\frac{d+2}{2(d+1)}\alpha^d + d+3 - {\alpha + d \choose d}}\\
\quad \quad  \geq \displaystyle{ \alpha^{d-1}\left[\frac{d+2}{2(d+1)}\alpha-\frac{(\alpha 
+d)(d+1)}{2d^2}\right]  - \frac{d+2}{d}\alpha + 1} \\ \quad  \displaystyle{
\quad =  \frac{\alpha (d^3+d^2 - 2d-1) - d(d+1)^2}{2d^2(d+1)}\alpha^{d-1} - 
\frac{d+2}{d}\alpha + 1}\\ \quad \displaystyle{
\quad \geq  \frac{d (d^3+d^2 - 2d-1) - d(d+1)^2}{2d^2(d+1)}\alpha^{d-1} - 
\frac{d+2}{d}\alpha + 1} \ (\text{since}\ \alpha \geq d\geq 3) \\ \quad \displaystyle{
\quad =  \frac{d^3- 4d-2 }{2d(d+1)}\alpha^{d-1} - \frac{d+2}{d}\alpha + 1} =: B
\end{array}$$
If $\alpha =3$, then $d=3$ and $B= 7/8$. For $\alpha \geq 4$, since $\alpha^{d-1}\geq 
4\alpha $, we further get
$$B \geq \frac{2(d^3- 4d-2 )}{d(d+1)}\alpha - \frac{d+2}{d}\alpha + 1= \frac{d[d(2d-1)-11]-6 
}{d(d+1)}\alpha + 1 >1.$$
Thus we always have $B>0$, which proves (\ref{EB3}).

Now we consider the case when an edge of $\Pcal$ is full, i.e. there are exactly $\alpha +1$ 
points on it belonging to $\Acal$. If $\alpha \geq 6$, then $\alpha \geq \frac{4d}{d-1}$ and 
the second condition is satisfield, so we are done. Since $\alpha \geq d$, the left cases are 
$d=4, \ \alpha \leq 5$ and $d=3,\ \alpha =4,5$. In these cases, by Lemma \ref{A2}(ii), $r(S) 
\leq \alpha^{d-2}$. 

If $d=4,\ \alpha \leq 5$, then $\deg K[S] - c+1 \geq \alpha^3 - {\alpha +3 \choose 3} + 6 > 
2\alpha^2\geq 2r(S)$, and by Theorem \ref{A1} we are done.

If $d=3,\ \alpha \leq 5$, let $\tilde{c} = \sharp(M_{\alpha,3}\setminus \Acal)$. Then  
$r(S)\leq \alpha $, and the inequality
$$ \deg K[S] - c+1 = \alpha^2 - {\alpha +2 \choose 2} + \tilde{c} + 4 \geq 2\alpha -1$$
does not hold only in the following situations: $\alpha =3,4,\ \tilde{c}=1,2$ and $\alpha =5, \ 
\tilde{c}=1$. By Theorem \ref{A1}, we can restrict ourselves to these situations. By Corollary 
\ref{C4}, we may assume that one deleting point is $(\alpha -1, 1,0)$. Thus, in each case there 
are only few configurations to consider. Using computer, we can check that $r(S)= 2$ if 
$\alpha =3,5$, and $r(S) \leq 3$ if $\alpha =4$. But then $\deg K[S] - c+1 = \alpha^2 - {\alpha 
+2 \choose 2} + \tilde{c} + 4 \geq 2r(S)-1$. Again by Theorem \ref{A1} we are done.
 \hfill $\square$
\end{pf}

{\bf Acknowledgment}.  All computations  in this paper were done by using the package 
CoCoa \cite{CCA}.

\end{document}